\documentclass[10pt]{article}
\usepackage{graphicx}
\usepackage{amsmath}
\usepackage[dvips]{epsfig}
\usepackage{amssymb}

\makeatletter
\renewcommand{\section}{\@startsection
  {section}%
  {2}%
  {0mm}%
  {\baselineskip}%
  {0.3 \baselineskip}%
  {\centering}}
\makeatother
\begin{document}

\title { A note on the alternating sums of powers of
consecutive $q$-integers }
\author{ Taekyun Kim$^{\dagger}$ \\[0.5cm]
         ${}$
${\dagger}$ Jangjeon Research Institute for Mathematical Science and Physics,\\
        544-4 Young-Chang Ri Hapcheon-Gun Kyungnam, 678-802, S. Korea\\
          {\it e-mail: tkim $@$kongju.ac.kr/ tkim64$@$hanmail.net }
        \\ \\ }

\date{}
\maketitle

{\footnotesize {\bf Abstract}\hspace{1mm}
  In this paper we construct a new $q$-Euler numbers and polynomials.
  By using these numbers and polynomials, we give the interesting
   formulae related to alternating sums of powers of consecutive $q$-integers
 following an idea due to Euler.


\bigskip
\section{Introduction }
\indent The Bernoulli numbers among the most interesting and
important number sequence in mathematics. They appeared in the
posthumous work `` ARS Conjectandi (1713)" by Jacob
Bernoulli(1654-1705) in connection with sums of powers of
consecutive integers(1713). Let $n, k$ be positive integers, and
let $S_{n,q}(k)$ be the sums of the $n$th powers of positive
integers up to $k-1$: $S_{n}(k)=\sum_{l=0}^{k-1} l^n $. Then it
was known that the sums of powers of consecutive integers due to
J. Bernoulli as follows:
$$S_n(k)=\frac{1}{n+1} \sum_{i=0}^n \binom {n+1}i B_ik^{n+1-i},
\text{ where $B_n$ are the $n$th Bernoulli numbers , cf.[7, 10,
15]}.$$ In [1], Carlitz has introduced an interesting $q$-analogue
of Bernoulli numbers, $\beta_{k, q}$. He has indicated a
corresponding Stadudt-Clausen theorem and also some interesting
congruence properties of the $q$-Bernoulli numbers. Let $q$ be an
indeterminate which can be considered in the complex number field,
and for any integer $k$ define the $q$-integer as
$[k]_q=\frac{1-q^k}{1-q}$, cf. [2, 5, 11]. Note that
$\lim_{q\rightarrow 1}[k]_q =k.$ For any positive integers $n, k$,
let $S_{n,q}(k)=\sum_{l=0}^{k-1}q^l[l]_q^n$, sf. [7, 10]. Then we
evaluated sums of powers of consecutive $q$-integers as follows:
$$S_{n,q}(k)=\frac{1}{n+1}\sum_{i=0}^n\binom{n+1}i
\beta_{i,q}q^{ki}[k]_q^{n+1-i}-\frac{(1-q^{(n+1)k})\beta_{n+1,q}}{n+1},
\text{ cf, [7, 10]}.$$

 The ordinary Euler numbers are defined by
$$
\dfrac{2}{e^{t}+1} = e^{Et} = \sum_{n=0}^{\infty} E_{n}
\dfrac{t^n}{n!}, \quad {\rm cf.}~ [2, 3, 5, 6, 7, 8],
$$
where we use the usual convention about replacing $E^{n}$ by
$E_{n} (n > 0)$ symbolically. Let $n, k$ be positive integers, and
let $T_{n,q}(k)$ be the alternating sums of the $n$th powers of
positive integers up to $k-1$: $T_{n}(k)=\sum_{l=0}^{k-1}
(-1)^ll^n $. Then Euler investigated the below formulas:
$$T_n(k)=\frac{(-1)^{n+1}}{2} \sum_{l=0}^{n-1} \binom {n}l E_lk^{n-1}+\frac{E_n}{2}\left(1+(-1)^{k+1}\right).$$

 Let $u$ be algebraic in complex
number field. Then Frobenius-Euler numbers are defined by
$$\dfrac{1-u}{e^{t}-u} =  \sum_{n=0}^{\infty} H_{n}(u)
\dfrac{t^n}{n!}, \quad {\rm cf.} ~[3, 5],
$$ note that $H_{n}(-1)=E_{n}$, cf. [2, 3, 5, 6, 7].
Carlitz has also introduced an interesting $q$-analogue of
Frobenius-Euler numbers in [1].  A recent author's study of more
general $q$-Euler numbers are found in previous publication [14].
In [4] we gave the new construction of $q$-Euler numbers,
$E_{n,q}^{*}$, which are different than Carlitz's $q$-extension
and author's $q$-extension in previous publication ( see [11]).
Let $\mathbb{Z}_{p}$ be the ring of $p$-adic integers, and let $p$
be a fixed odd prime number. Then the $p$-adic $q$-integral was
defined by author as follows:
$$\int_{\mathbb{Z}_{p}} f(x) d\mu_{q} (x) = \lim_{N \rightarrow
\infty} \dfrac{1}{[p^{N}]_{q}} \sum_{x=0}^{p^{N}-1} f(x)q,$$ where
$f \in UD(\mathbb{Z}_{p})$, cf. [5, 8, 9, 10, 11]. The above
$q$-extension of Euler numbers, $E_{n,q}^{*}$, were written by
\begin{eqnarray}
\sum_{n=0}^{\infty} E_{n,q}^{*} \dfrac{t^n}{n !} &=&
\sum_{n=0}^{\infty} \int_{\mathbb{Z}_{p}} [x]_{q}^{n} d\mu_{-q}(x)\dfrac{t^{n}}{n!}\notag\\
&=& [2]_{q} \sum_{n=0}^{\infty} (-1)^{n} q^{n} e^{[n]_{q}t}, \quad
({\rm see} ~[4, 6]).
\end{eqnarray}
By (1), we easily see that
$$E_{n, q}^{*} = [2]_{q} \Big(\dfrac{1}{1-q} \Big)^{n}
\sum_{l=0}^{n} \binom{n}{l} (-1)^{l} \dfrac{1}{1+q^{l+1}}, \quad
({\rm see} ~[4, 6]). $$
 For $n, m \in\Bbb N$, we gave the below interesting formula:
$$\sum_{l=0}^{n-1}(-1)^lq^l[l]_q^m=\frac{1}{[2]_q}\left((-1)^{n+1}q^nE_{m,q}^*(n)+E_{m,q}^*\right), \text{ see [4]}.$$
In this paper we consider a new approach to $q$-Euler numbers and
polynomials and give some identities and properties between
$q$-Euler numbers and polynomials. Finally we will evaluate the
value of $\sum_{l=0}^{n-1}(-1)^l[l]_q^m$ by using our new
$q$-Euler numbers and polynomials. This formula seems to be nice.

\bigskip

\section{ A note on q-Euler numbers and polynomials }

For $q \in \mathbb{C}$ with $|q| < 1$, we consider a modified
$q$-extension of Euler numbers as follows:
\begin{equation}
2\sum_{l=0}^{\infty} (-1)^{l} e^{{[l]_{q}} t} =
\sum_{n=0}^{\infty} E_{n, q}\dfrac{t^{n}}{n!}.
\end{equation}
From (2), we can derive the below formula ;
\begin{eqnarray}
2 \sum_{l=0}^{\infty} (-1)^{l} e^{{[l]_q}t} &=& 2
e^{\dfrac{t}{1-q}} \sum_{l=0}^{\infty} (-1)^{l}
\sum_{j=0}^{\infty} \big(\dfrac{1}{1-q}\big)^{j} q^{lj}
\dfrac{t^{j}}{j!}\notag\\
&=& 2 e^{\dfrac{t}{1-q}} \sum_{j=0}^{\infty}
\big(\dfrac{1}{1-q}\big)^{j} (-1)^{j}
\dfrac{1}{1+q^{j}}\dfrac{t^{j}}{j!}\notag\\
&=& 2 \sum_{i=0}^{\infty} \big(\dfrac{1}{1-q}\big)^{i}
\dfrac{t^{i}}{i!} \sum_{j=0}^{\infty} \big(\dfrac{1}{1-q}\big)^{j}
(-1)^{j} \dfrac{1}{1+q^{j}} \dfrac{t^{j}}{j!}\notag\\
&=& 2 \sum_{\substack{n=0 \\ n = i+j}}^{\infty} \Big(
\big(\dfrac{1}{1-q}\big)^{n} \sum_{j=0}^{n} (-1)^{j}
\big(\dfrac{1}{1+q^{j}}\big) \dfrac{n!}{j! (n-j)!}\Big)
\dfrac{t^{n}}{n!}\notag\\
&=& \sum_{n=0}^{\infty}\Big(2\big(\dfrac{1}{1-q}\big)^{n}
\sum_{j=0}^{n} \binom{n}{j} (-1)^{j} \dfrac{1}{1+q^{j}}
\Big)\dfrac{t^{n}}{n!}.
\end{eqnarray}
Thus we have the following :
\vspace{0.1in}\\
\noindent {\bf Theorem 1.} {\it For $n \geq 0$, we have $$E_{n, q}
= 2 \big(\dfrac{1}{1-q} \big)^{n} \sum^{n}_{j=0} \binom{n}{j}
(-1)^{j} \dfrac{1}{1+q^{j}}.$$ Note that $\underset{q\rightarrow
1}{\lim} E_{n, q} = E_{n}$.}

\medskip

By simple calculation, it is easy to check that
\begin{eqnarray}
2 e ^{t[x]_{q} }\sum_{l=0}^{\infty} (-1)^{l} e^{[l]_{q}q^{x} t}
&=& 2 \sum_{l=0}^{\infty} (-1)^{l} e^{([x]_{q}+[l]_{q} q^{x}) t}\notag\\
&=& 2 \sum_{l=0}^{\infty} (-1)^{l} e^{[x+l]_{q} t}.
\end{eqnarray}

From (4), we can define the below $q$-Euler polynomials :
\begin{equation}
2 \sum_{l=0}^{\infty} (-1)^{l} e^{[x+l]_{q} t}=
\sum_{n=0}^{\infty} E_{n, q}(x) \dfrac{t^{n}}{n!}.
\end{equation}

By (4) and (5), we easily see that
\begin{eqnarray}
\sum_{n=0}^{\infty} E_{n, q}(x) \dfrac{t^{n}}{n!} &=& 2
\sum_{l=0}^{\infty} (-1)^{l} \sum_{j=0}^{\infty} [x+l]_{q}^{j}
\dfrac{t^{j}}{j!}\notag\\
&=& 2 \sum_{l=0}^{\infty} (-1)^{l} \sum_{n=0}^{\infty}
\binom{n}{j}
(-1)^{j} q^{xj} q^{lj} \dfrac{t^{n}}{n!}\notag\\
&=& 2\sum_{n=0}^{\infty} \big(\dfrac{1}{1-q})^{n} \sum_{j=0}^{n}
\binom{n}{j} (-1)^{j} q^{xj} \dfrac{1}{1+q^{j}}\dfrac{t^{n}}{n!}.
\end{eqnarray}

From (4), (5) and (6), we can obtain the following :
\vspace{0.1in}\\
\noindent {\bf Theorem 2.} {\it For $n \geq 0$, we have
\begin{eqnarray*}
E_{n, q}(x) &=& 2 \big(\dfrac{1}{1-q} \big)^{m} \sum^{n}_{j=0}
\binom{n}{j} (-1)^{j} q^{xj} \dfrac{1}{1+q^{j}} \\
&=&
\sum_{k=0}^{\infty} \binom{n}{k} q^{kx} E_{k, q} [x]_{q}^{n-k}.
\end{eqnarray*}}
By (2), (4) and (5), we easily see that
\begin{eqnarray*}
&&\sum_{m=0}^{\infty} \big( (-1)^{m+1} E_{m, q}(n) + E_{m, q}
\big) \dfrac{t^{m}}{m!}\\
&=& -2 \sum_{l=0}^{\infty} (-1)^{l+n} e^{[l+n]_{q} t} + 2
\sum_{l=0}^{\infty} (-1)^{l} e^{[l]_{q} t}\\
&=& 2 \sum_{l=0}^{n-1} (-1)^{l} e^{[l]_{q} t}.
\end{eqnarray*}
Therefore we obtain the following theorem:
\vspace{0.1in}\\
\noindent {\bf Theorem 3.} {\it Let $m$ be the positive integers
bigger than 1. Then we have
\begin{equation*}
\hspace{1in}\dfrac{\big((-1)^{1+m} E_{m, q}(n) + E_{m, q}
\big)}{2} = \sum_{l=0}^{n-1} (-1)^{l}
[l]_{q}^{m}.\hspace{1in}(6.1)
\end{equation*}}
\vspace{0.1in}\\
\noindent {\bf Remark.} In [4], it was known that
$$
\hspace{0.85in}\dfrac{1}{[2]_{q}} \big((-1)^{1+m} E_{m, q}^{*}
(n)+E_{m, q}^{*}\big)= \sum_{l=0}^{n-1} (-1)^{l} q^{l}
[l]_{q}^{m}.\hspace{0.85in}(6.2)$$

Comparing (6.1) and (6.2), we see that one new formula (6.1) also
seems worthwhile and valuable as same as (6.2)
\vspace{0.1in}\\
\noindent {\bf Remark.} In [5], we note that
\begin{eqnarray*}
\sum_{n=0}^{\infty} E_{n} (x) \dfrac{t^{n}}{n!} &=&
\dfrac{2}{e^{t}+1} e^{xt}\\
&=& \lim_{q\rightarrow 1} \big(2\sum_{l=0}^{\infty} (-1)^{l}
e^{[x+l]_{q} t}\big) \\
&=& \sum^{\infty}_{n=0} \big(\lim_{q \rightarrow 1} E_{n, q}(x)
\big) \dfrac{t^{n}}{n!}.
\end{eqnarray*}

Hence, we have $\underset{q\rightarrow 1}{\lim} E_{n, q}(x) =
E_{n}(x)$, where $E_{n}(x)$ are called ordinary Euler polynomials.

From(5), we can also derive ($f \in \mathbb{N}$, odd)
\begin{eqnarray*}
\sum_{n=0}^{\infty} E_{n, q}(x) \dfrac{t^{n}}{n!} &=& 2
\sum_{n=0}^{\infty} (-1)^{n} e^{[x+n]_{q} t}\\
&=&  2
\sum_{n=0}^{\infty} \sum_{a=0}^{f-1} (-1)^{a+nf} e^{[x+a+nf]_{q} t}\\
&=&  \sum_{a=0}^{f-1} (-1)^{a} \Big(2 \sum_{n=0}^{\infty} (-1)^{n} e^{[f]_{q} [\frac{x+a}{f}+n]_{qf} t}\Big)\\
&=&  \sum_{m=0}^{\infty} \Big( \sum_{a=0}^{f-1} (-1)^{a} E_{m,
q^f} \big(\dfrac{x+a}{f}\big)\Big) [f]^{m}_{q} \dfrac{t^{m}}{m!},
\end{eqnarray*}
where $f$ is odd positive integer.

By comparing the coefficients on both sides , we obtain the
following theorem.
\vspace{0.1in}\\
\noindent {\bf Theorem 4.} {\it Let $f$ be a positive odd integer.
Then we have
$$[f]_{q}^{n} \sum_{a=0}^{f-1} (-1)^{a} E_{m, q^{f}}
\big(\dfrac{x+a}{f} \big) = E_{m, q}(x).$$}

For $s\in \mathbb{C}$, let us consider the following complex
integration :
\begin{equation} \dfrac{1}{\Gamma(s)}
\int_{0}^{\infty} \sum_{l=0}^{\infty} (-1)^{l} e^{-[x+l]_{q} t} dt
= \sum_{l=0}^{\infty} \dfrac{(-1)^{l}}{[x+l]_{q}^{s}}.
\end{equation}
Thus, we can define the Euler $q$-zeta function as follows :
\begin{equation}
\zeta_{E, q} (s, x) =\sum_{n=0}^{\infty}
\dfrac{(-1)^{n}}{[n+x]_{q}^{s}}, \text{ $s\in \mathbb{C}$}.
\end{equation}

By (5), (7) and (8), we easily see that $\zeta_{E, q}(-n, x) =
\frac{1}{2} E_{n, q}(x), \quad n \in \mathbb{N}$.
\vspace{0.1in}\\
\noindent {\bf Remark.} Let $E_{n}(x)$ be the ordinary Euler
polynomials . Then we know that
\begin{equation}
\dfrac{\big( (-1)^{m+1} E_{m}(n) +E_{m} \big)}{2} =
\sum_{l=0}^{n-1} (-1)^{l} l^{m}, \quad {\rm see} ~[12].
\end{equation}
Theorem 3 is the new $q$-extension of Eq.(9). Let
$$H_q(s, a; F)=\sum_{m\equiv a (F), m>0}\frac{(-1)^m}{[m]_q^s}=\sum_{n=0}^{\infty}
\frac{(-1)^{a+nF}}{[a+nF]_q^s}=[F]_q^{-s}(-1)^a\zeta_{E,q^F}(s,
\frac{a}{F}),$$ where $a$ and $F$(=odd) are positive inteegrs with
$0<a<F$. Then we have
$$H_q(-n, a;
F)=\frac{(-1)^a[F]_q^nE_{n,q^F}(\frac{a}{F})}{2}, \text{ $n \geq
1$.}$$ Let $\chi$ be the Dirichlet character with conductor $d\in
\Bbb N$(=odd). Then we define the generalized $q$-Euler numbers
attached to $\chi$ as follows:
$$F_{\chi,q}(t)=2\sum_{n=0}^{\infty}e^{[n]_qt}\chi(n)(-1)^n=\sum_{n=0}^{\infty}
E_{n,\chi,q}\frac{t^n}{n!}.$$ Note that
$$E_{n,\chi, q}=[d]_q^n\sum_{a=0}^{d-1}\chi(a)(-1)^a E_{n,
q^d}(\frac{a}{d}).$$ For $s\in\Bbb C,$ let us define the
$q-l-$function as follows:
$$l_{E, q}(s,
\chi)=\sum_{n=1}^{\infty}\frac{(-1)^n\chi(n)}{[n]_q^s}=\frac{1}{2}\frac{1}{\Gamma
(s)}\int_{0}^{\infty}F_{\chi, q}(-t)t^{s-1}dt. $$ Then we easily
see that $l_{E,q}(-n, \chi)=\frac{1}{2}E_{n,\chi,q},$ ( $n\in\Bbb
N$). For $\in\Bbb C ,$ it is easy to see that
$$l_{E,q}(s,\chi)=\sum_{a=1}^F\chi(a)H_q(s,a; F). $$
The function $H_q(s,a;F)$ will be called the partial Euler
$q$-zeta function. Finally we suggest the below problem.

\noindent {\bf Problem.} Find the Witt's formula for the $q$-Euler
numbers $(E_{n}, q)$, which was defined in this paper. In [5], the
Witt's type formula for $E_{n, q}^{*}$ was given by
\begin{eqnarray*}
\sum_{m=0}^{\infty} E_{n, q}^{*} \dfrac{t^{n}}{n!} &=&
\sum_{n=0}^{\infty} \int_{\mathbb{Z}_{p}} [x]_{q}^{n} d\mu_{-q}(x)
\dfrac{t^{n}}{n!}\\
&=& [2]_{q} \sum_{n=0}^{\infty} (-1)^{n} q^{n} e^{[n]_{q}t}.
\end{eqnarray*}
By the same method, it seems to be possible that we give the
Witt's formula of $q$-Euler numbers which can be represented by
$p$-adic $q$-integrals as follows:
\begin{eqnarray*}
\sum_{m=0}^{\infty} E_{n, q} \dfrac{t^{n}}{n!} &=&
\sum_{n=0}^{\infty} \int_{\mathbb{Z}_{p}} [x]_{q}^{n}
d\mu_{\Box}(x)\\
&=& 2 \sum_{l=0}^{\infty} (-1)^{l}  e^{[l]_{q} t}, \quad ({\rm
cf.} ~[4, 6, 14]).
\end{eqnarray*}
\vspace{0.2in}\\
\footnotesize{

\end{document}